\documentclass[11pt]{article}

\usepackage{latexsym, amssymb}

\newcommand{\be}{\begin{equation}}
\newcommand{\ee}{\end{equation}}
\newcommand{\beq}{\begin{eqnarray}}
\newcommand{\eeq}{\end{eqnarray}}

\newtheorem{prop}{Proposition}

\def\bG{\bar{G}}
\def\hG{\hat{G}}
\def\hg{\hat{g}}
\def\tg{\tilde{g}}
\def\tM{\tilde{M}}
\def\stop{\hfill$\Box$}

\def\S{\Sigma}
\def\pf{\noindent {\em{Proof}}: }

\def\R{\mathbb{R}}

\def\bg{\bar{g}}
\def\mS{\mathbb{S}^3}
\def\ep{\epsilon}

\def\RP{{RP}}

\def\vh{\vspace{.3cm}}

\begin{document}

\title{A Note on Existence and Non-existence of Minimal Surfaces in 
Some Asymptotically Flat $3$-manifolds}

\author{Pengzi Miao\thanks{%
Current address: Department of Mathematics, University of California, Santa Barbara, CA 93106, USA.  E-mail: pengzim@math.ucsb.edu}
\thanks{ Address after April 3 06: School of Mathematical Sciences, 
Monash University,  Victoria 3800, Australia.} 
 }

\date{}

\maketitle

\begin{abstract}
%We consider conformal classes $[g]$ of metrics on $S^3$ which have a metric of positive Ricci curvature. Defining $V_{max}$ to be the supremum of the volume of $\bg \in [g]$ satisfying $Ric(\bg) \geq 2 \bg$, we show if  $V_{max}$ is no less than one half of the volume of the standard unit sphere, then the asymptotically flat manifold $(S^3 \setminus \{P\}, G^4 g)$ has no closed minimal surfaces. Here $G$ is the  Green's function of the conformal Laplacian of $(S^3, g)$ at a point $P$.  We also give an example of $[g]$ where $(S^3 \setminus \{P\}, G^4 g)$ does have closed minimal surfaces.

Motivated by problems on apparent horizons in general relativity, we prove the following theorem on minimal surfaces: Let $g$ be a metric on the three-sphere $S^3$ satisfying $Ric(g) \geq 2 g$. If the volume of $(S^3, g)$ is no less than one half of the volume of the standard unit sphere, then there are no closed minimal surfaces in the asymptotically flat manifold $(S^3 \setminus \{ P \}, G^4 g)$.  Here $G$ is the  Green's function of the conformal Laplacian of $(S^3, g)$ at an arbitrary point $P$. We also give an example of $(S^3, g)$ with $Ric(g) > 0$ where $(S^3 \setminus \{ P \}, G^4 g)$ does have closed minimal surfaces.

%We consider asymptotically flat manifolds of the form $(S^3 \setminus \{ P \}, G^4 g)$, where $G$ is the Green's function of the conformal Laplacian of $(S^3, g)$ at a point $P$. We show  if  $Ric(g) \geq 2 g$ and the volume of $(S^3, g)$ is no less than one half of the volume of the standard unit sphere, then there are no closed minimal surfaces in $(S^3 \setminus \{ P \}, G^4 g)$. We also give an example of $(S^3, g)$ where $Ric(g) > 0$ but $(S^3 \setminus \{ P \}, G^4 g)$ does have closed minimal surfaces. 

\end{abstract}

\section{Introduction}
Let $(N^3, g, p)$ be an initial data set satisfying the dominant energy constraint condition in general relativity. It is a fascinating question to ask under what conditions 
an {\em apparent horizon} (of a back hole) exists in $(N^3, g, p)$. Here an apparent horizon is a $2$-surface $\S^2 \subset N^3$ satisfying
\be \label{eqofhorizon}
H_\S = \mathrm{Tr}_\S p, 
\ee
where $H_\S$ is the mean curvature of $\S$ in
$N$ and $\mathrm{Tr}_\S p$ is the trace of the restriction of $p$
to $\S$. 

A fundamental result of Schoen and Yau states  that
 {\em matter condensation} causes apparent horizons to be formed 
\cite{Schoen-Yau-matter}. Their result is remarkable not only because it provides a general  criteria to the existence question, but also because it leads to a refined problem -- besides matter fields, what is the {pure effect of gravity} on the formation of apparent horizons? 

To analyze this refined problem, one considers an {asymptotically flat} initial data set $(N^3, g, p)$ in a {\em vacuum} spacetime. 
As the first step, one assumes $(N^3, g, p)$ is  {time-symmetric} (i.e. $p \equiv 0$).   In this context, an apparent horizon is  
simply a {\em minimal surface},  and the relevant topological assumption is 
that $N^3$ is diffeomorphic to $\R^3$. (If $N^3$ has nontrivial topology, 
a closed minimal surface always exists by \cite{Meeks-Simon-Yau}.) 

There is a geometric construction of such an initial data set.  
Let $[g]$ be a conformal class of metrics on the three-sphere $S^3$. 
Recall the Yamabe constant of $(S^3, [g])$ is defined by
\be
Y(S^3, [g]) = \inf_{ v \in W^{1,2}(S^3) }  
\frac{ \int_M [ 8 | \nabla v|^2_g +  R(g) v^2 ] dV_{g} }{
\left( \int_M v^6 dV_{g} \right)^\frac{1}{3} },
\ee
where $R(g)$ is the scalar curvature of $g$. If $Y(S^3, [g]) > 0$, 
there exists a positive Green's function $G$ of the conformal Laplacian
$8\triangle_g -  R(g) $ at any fixed point $P  \in S^3$.
Consider the new metric $G^4 g$ on $S^3 \setminus \{P\}$, it is easily
checked that $(S^3 \setminus \{ P \}, G^4 g)$ is {asymptotically flat}
with {zero scalar curvature}. One basic fact about this construction
is that the blowing-up manifold $(S^3 \setminus \{ P\}, G^4g)$,
up to a constant scaling, depends only on the {conformal class} $[g]$. 
Precisely, if one replaces $g$ by another metric $\bg \in [g]$ and
let $\bar{G}$ be the Green's function associated to $\bg$, then the metric $\bar{G}^4 \bg$  differs from  $ G^4g$ only by a constant multiple. Therefore, it is of interest to seek conditions on $[g]$ 
that determine whether $(S^3 \setminus \{ P\}, G^4 g)$ has a horizon. 

So far, no such a conformal invariant condition has been found. However, there are results where conditions in terms of a single metric are given. In \cite{Beig-OMurchadha}, 
Beig and \'{O} Murchadha studied the behavior of  a {\em critical sequence}, i.e. a sequence of metrics $\{g_n\}$ on $S^3$ converging to a metric $g_0$ with zero scalar curvature. They showed 
the blowing-up manifold  $(S^3 \setminus \{P\}, G_n^4 g_n)$ has a horizon for sufficiently large $n$. Their idea was further explored by Yan \cite{Yan}. Given a metric $g$ on $S^3$, assuming the diameter of $(S^3, g) \leq D$, the volume of $(S^3, g) \geq V$ and the Ricci curvature of $g$ satisfies $Ric(g) \geq \mu g$, Yan showed
that, for any $r > \frac{3}{2}$, there exists a small positive number 
$\delta=\delta(\mu, V, D, r)\leq1$ such that,  if
$R(g) > 0$ and $||R(g)||_{L^r(S^3, g)} < \delta$, then the blowing-up 
manifold $(S^3 \setminus \{P\}, G^4 g)$ has a horizon. 

One  question arising from Yan's theorem is whether a {\em
positive} Ricci curvature metric on $S^3$ can produce a blowing-up 
manifold  with a horizon, as it is unclear whether Yan's theorem could be applied when $\mu > 0$. Another motivation to this question is,
as a positive Ricci curvature metric can be deformed to the standard metric on $S^3$ through metrics of positive Ricci curvature,
it is of potential interest to study how the horizon disappears in the 
corresponding deformation of the blowing-up manifold if it
exists initially. 

In this paper, we focus on conformal classes of metrics with a positive Ricci curvature metric. Our main result is the observation of a  volume
condition which guarantees non-existence of horizons in the blowing-up manifold. Throughout the paper, $\mS$ denotes $S^3$ with
the standard metric of constant curvature $+1$.

\vh

\noindent {\bf Theorem}
{\em Let $[g]$ be a conformal class of metrics on $S^3$ which 
has a metric of positive Ricci curvature. Consider 
$$
V_{max}(S^3, [g]) = 
\sup_{ \bg \in [g]} \{ Vol(S^3, \bg) \ | \ Ric(\bg) \geq 2 \bg \} ,
$$
where $Vol(\cdot)$ is the volume functional.
If 
$$ 
V_{max}(S^3, [g]) \geq \frac{1}{2} Vol(\mS), 
$$
then the asymptotically flat manifold $(S^3 \setminus \{P\}, G^4 g)$
has no horizon.}

\vh

We also give an example of $(S^3, g)$ with $Ric(g)>0$ where
$(S^3 \setminus \{P\}, G^4 g)$ does have horizons.

\section{Positive Ricci curvature and maximum volume}
We first explain the volume assumption in the Theorem.
Let $M^n$ be a smooth, connected, closed manifold of dimension $n \geq 3$. 
Assume $[g]$ is a conformal class of metrics on $M^n$ which has
a metric of positive Ricci curvature.  One can define
\be
V_{max}(M^n, [g]) = \sup_{ \bg \in [g]} \{ Vol(M^n, \bg) \ | \ 
Ric(\bg) \geq (n-1) \bg \} .
\ee
The following result relating $V_{max}(M^n, [g])$ and the Yamabe 
constant of $(M^n, [g])$ was observed in \cite{Guan-Wang}.

\begin{prop}
Let $[g]$ be a conformal class of metrics on $M^n$ which has 
a metric of positive Ricci curvature.  Then the Yamabe constant of 
$(M^n, [g])$ satisfies
\be
Y(M^n, [g]) \geq n(n-1) V_{max}(M^n, [g])^\frac{2}{n} .
\ee
\end{prop}

\pf By definition, 
\be
Y(M^n, [g]) = \inf_{ v \in W^{1,2}(M) }  
\frac{ \int_M [ c_n | \nabla v|^2_{\bg} +  R(\bg) v^2 ] dV_{\bg} }{
\left( \int_M v^\frac{2n}{n-2} dV_{\bg} \right)^\frac{n-2}{n} } 
\ee
for any $\bg \in [g]$, where $c_n = \frac{4(n-1)}{n-2}$.

Now we assume $Ric(\bg) \geq (n-1) \bg$. Then by a result of Ilias \cite{Ilias}, which is based on the isoperimetric inequality of Gromov \cite{Gromov-isoineq}, we have
\be
 \int_M [ c_n | \nabla v|^2_{\bg}  + n(n-1) v^2 ] dV_{\bg}
\geq
\left( \int_M v^\frac{2n}{n-2} dV_{\bg} \right)^\frac{n-2}{n}
n(n-1) Vol(M^n, \bg)^\frac{2}{n} 
\ee
for any $v \in W^{1,2}(M)$. Note that $R(\bg) \geq n(n-1)$, hence
\beq
Y(M^n, [g]) & \geq & \inf_{ v \in W^{1,2}(M) }  
\frac{ \int_M [ c_n | \nabla v|^2_{\bg} +  n (n-1) v^2 ] dV_{\bg} }{
\left( \int_M v^\frac{2n}{n-2} dV_{\bg} \right)^\frac{n-2}{n} }
\nonumber \\
& \geq & n(n-1) Vol(M^n, \bg)^\frac{2}{n} .
\eeq
Taking the supremum over $\bg \in [g]$ satisfying $Ric(\bg) \geq (n-1) \bg$,  
we have
\be
Y(M^n, [g]) \geq n(n-1) V_{max}(M^n, [g])^\frac{2}{n} .
\ee 
\stop

\vh

As an immediate corollary, we see the assumption 
$$ V_{max}(S^3, [g]) \geq \frac{1}{2} Vol(\mS) $$
in the Theorem implies 
\beq \label{volYamabe}
Y(S^3, [g]) & \geq & 6 \left( \frac{1}{2} \right)^\frac{2}{3} Vol(\mS)^\frac{2}{3}
\nonumber \\
& = & Y(RP^3, [g_0]),
\eeq
where $RP^3$ is the three dimensional projective space and $g_0$
is the standard metric on $RP^3$ which has constant sectional curvature $+1$.

\section{An upper bound of the Sobolev constant when a horizon is present}
One basic fact relating the conformal class $[g]$ on $S^3$ and the blowing-up metric $ h =G^4 g$ on $\R^3 = S^3 \setminus \{ P\}$ is 
\be \label{Yamabeandsob}
Y(S^3, [g]) = 8 S( h ),
\ee
where $S( h )$ is the Sobolev constant of the asymptotically flat manifold $(\R^3, h)$ \cite{Bray-Neves}. 
Recall $S( h )$ is defined by
\be \label{defofsob}
S( h ) = \inf_{u \in W^{1,2}(\R^3, h)} \left\{ 
\frac{ \int_{\R^3} | \nabla u|^2_h \ dV_h }
{ (\int_{\R^3} u^6 \ d V_h )^\frac{1}{3}} 
\right\} .
\ee 
The next proposition, which plays a key role in the derivation of the
Theorem, was essentially established by Bray and Neves in \cite{Bray-Neves} using the inverse mean curvature flow technique \cite{IMF}.
As the statement of Bray and Neves is  different from what we 
need, we include the proof here.

\begin{prop} \label{sob}
Let $h$ be a complete metric on $\R^3$ such that $(\R^3, h)$ is 
asymptotically flat. If $(\R^3, h)$ has nonnegative scalar curvature
and has a closed minimal surface,
then 
\be \label{estofsob}
S (h)< \frac{1}{8} Y(RP^3, [g_0]).
\ee 
\end{prop}

\pf 
Since $(\R^3, h)$ has a closed minimal surface, 
the {\em outermost} minimal surface $\mathcal{S}$ in 
$(\R^3, h)$, i.e. the closed minimal surface that is not enclosed by 
any other minimal surface \cite{Bray_Penrose},  exists and consists 
of a finite union of disjoint, embedded minimal two-spheres and 
projective planes. As our background manifold is $\R^3$, 
$\mathcal{S}$ must consist of embedded minimal two-spheres alone, 
furthermore each component of $\mathcal{S}$ 
necessarily bounds a three-ball.

We fix a component $\S$ of $\mathcal{S}$ and denote by $\Omega$ the three-ball that $\S$ bounds in $\R^3$.  Let $\phi$ be the weak solution to the {inverse mean curvature flow}  in $(\R^3 \setminus \bar{\Omega}, \ h)$ with initial condition $\S$ \cite{IMF}. $\phi$ satisefies
$$ \phi \geq 0, \ \phi|_\S = 0, \ \ \lim_{x \rightarrow \infty}\phi =
\infty .$$ 
Let $\S_t$ be the set $\partial \{ u < t \}$ 
for $t>0$ and $\S_0$ be the starting surface $\S$, 
then the family of surfaces $\{ \S_t \}$ satisfies 
the following properties \cite{IMF}:

\begin{enumerate}

\item $\{ \S_t \}$ consists of $C^{1, \alpha}$ surfaces. For a.e. $t$, 
$\S_t$ has weak mean curvature $H$ and 
$H = |\nabla u|_h$ for a.e. $x \in \S_t$.

\item $|\S_t| = e^t |\S_0|$, where $|\S_t|$ denotes the area of 
$\S_t$. 

\item Since $(\R^3, h)$ has  {nonnegative}
scalar curvature, $\S$ is {connected} and 
$\R^3 \setminus \bar{\Omega}$ is {simply connected}, 
the Hawking quasi-local mass of $\S_t$,
$$ m_H(\S_t) = \sqrt{ \frac{|\S_t|}{16 \pi} } \left( 
1 - \frac{1}{16 \pi} \int_{\S_t} H^2 d \mu
\right) ,$$
is monotone increasing. Here $d\mu$ is the induced surface measure. 
\end{enumerate} 
Now we restrict attention to functions $u \in W^{1,2}(\R^3, h)$ 
that have the form
\be
u (x) = \left\{
\begin{array}{cl}
f(0) & x \in \Omega \\
f(\phi(x)) & x \in \R^3 \setminus \Omega
\end{array}
\right.
\ee
for some $C^1$ functions $f(t)$ defined 
on $[0, \infty)$.
By the coarea formula and Property 1 above, we have
\beq
\int_{\R^3} | \nabla u|^2_h d V_h &  = &  \int_0^\infty f^\prime(t)^2 
\left( \int_{\S_t} H d \mu \right) dt \nonumber \\
& \leq & \int_0^\infty f^\prime(t)^2 
\sqrt{16 \pi |\S| (e^t - e^\frac{t}{2})}  dt ,
\eeq
where the inequality follows from Property 2, 3 and H\"{o}lder's inequality. Similarly, we have
\beq
\int_{\R^3} u^6 d V_h & \geq & \int_0^\infty f(t)^6 \left(
\int_{\S_t} H^{-1} d \mu \right)  dt \nonumber \\
& \geq & 
\int_0^\infty f(t)^6 e^{2t} | \S|^2
[16 \pi |\S| (e^t - e^\frac{t}{2})]^{-\frac{1}{2}} dt .
\eeq 
Therefore, 
\be \label{estofu}
\frac{\int_{\R^3} | \nabla u|^2_h d V_h }{
\left(\int_{\R^3} u^6 d V_h \right)^\frac{1}{3}} 
\leq 
\frac{(16 \pi)^\frac{2}{3} \int_0^\infty f^\prime(t)^2 
(e^t - e^\frac{t}{2})^\frac{1}{2}  dt }{
\left( \int_0^\infty f(t)^6 e^{2t} 
(e^t - e^\frac{t}{2})^{-\frac{1}{2}} dt 
\right)^\frac{1}{3} } .
\ee
To pick an optimal $f(t)$ that minimizes the right side
of (\ref{estofu}), we consider the half spatial Schwarzschild manifold
$$ (M^3, g_S) = 
( \R^3 \setminus B_1(0), ( 1 + \frac{1}{|x|} )^4 \delta_{ij} ) $$
and the quotient manifold $(\tM^3, \tg_S)$ obtained from 
$(M^3, g_S)$ by identifying the antipodal points of $\{ |x| = 1 \}$.
Up to scaling, $(\tilde{M}^3, \tilde{g}_S)$ is
isometric to $(\RP^3 \setminus \{ Q \}, G_0^4 g_0)$, 
the blowing-up manifold of $(\RP^3, g_0)$ by its Green function at 
a point $Q$. Hence, the Sobolev constant $S(\tg_S)$ of 
$(\tilde{M}^3, \tg_S)$ equals $\frac{1}{8} Y(\RP^3, [g_0])$. On the other hand,  $S(\tilde{g}_S)$ is achieved by a function  ${u}_0$ that is  constant  on each coordinate sphere $\{ |x| = t \}$ in $\tM$,
and the level set of the solution $\phi_0$ to the inverse mean curvature flow  starting at $\{ |x| = 1\}$ in $(M, g_S)$ is also given  by coordinate spheres.  Therefore, lifted as a function on $(M^3, g_S)$, $u_0$ has the form of 
$$ u_0 = f_0 \circ \phi_0  $$ for some explicitly determined function $f_0(t)$, and
\be
S(\tg_S) = \frac{\int_{M} | \nabla u_0|^2_{g_S} dV_{g_S} }
{ (\int_M u^6_0 \ dV_{g_S} )^\frac{1}{3} } = 
\frac{(16 \pi)^\frac{2}{3} \int_0^\infty f_0^\prime(t)^2 
(e^t - e^\frac{t}{2})^\frac{1}{2}  dt }{
\left( \int_0^\infty f_0(t)^6 e^{2t} 
(e^t - e^\frac{t}{2})^{-\frac{1}{2}} dt 
\right)^\frac{1}{3} } ,
\ee
where the second equality holds because the Hawking quasi-local
mass remains unchanged along the level sets of $\phi_0$.
Now consider $u = f_0 \circ \phi $ on $(\R^3, h)$. It was verified 
in \cite{Bray-Neves} that $u \in W^{1,2}(\R^3, h)$. Therefore, we have
\beq \label{estofsobineq}
S(h) \leq \frac{ \int_{\R^3} | \nabla u|^2_h \ dV_h }
{ (\int_{\R^3} u^6 \ d V_h )^\frac{1}{3}} 
& \leq &  \frac{(16 \pi)^\frac{2}{3} \int_0^\infty f_0^\prime(t)^2 
(e^t - e^\frac{t}{2})^\frac{1}{2}  dt }{
\left( \int_0^\infty f_0(t)^6 e^{2t} 
(e^t - e^\frac{t}{2})^{-\frac{1}{2}} dt 
\right)^\frac{1}{3} } \nonumber \\
& = &  S(\tg_S) = \frac{1}{8} 
Y(\RP^3, [g_0]).
\eeq

To show the strict inequality, we assume
 $S(h) =  \frac{1}{8}Y(\RP^3, [g_0])$. Then,
 $S(h)$ is achieved by $u=f_0 \circ \phi$. 
It follows from the Euler-Lagrange
equation of the Sobolev functional (\ref{defofsob}) that
$u$ satisfies
\be
\triangle_h u + C u^5 = 0 \ \ \mathrm{on} \ \R^3 ,
\ee
where $C = S(h) || u ||_{L^6(\R^3, h)}^{-4}$.
However, $u  \equiv f_0(0) $ on $\Omega$ and
$f_0(0) \neq 0$ (Indeed, up to a constant multiple, 
$ f_0(t) = ( 2 e^t - e^\frac{t}{2})^{-\frac{1}{2}} $ 
\cite{Bray-Neves}).
Hence, $C = 0$, which contradicts to the fact that 
$u$ is not a constant. Therefore,
the strict inequality $S(h) < \frac{1}{8}Y(\RP^3, [g_0])$ holds.  \stop

\vh

%The Theorem follows directly from (\ref{volYamabe}),
%(\ref{Yamabeandsob}) and Proposition \ref{sob}.

\noindent {\em Proof of the Theorem: }
Suppose $(S^3 \setminus \{P\}, G^4 g)$
has a horizon, then it follows from (\ref{Yamabeandsob})
and Proposition \ref{sob}  that
\be
 Y(S^3, [g]) <  Y (\RP^3, [g_0]) .
\ee
On the other hand, the assumption
$V_{max}(S^3, [g]) \geq \frac{1}{2} Vol(\mathbb{S}^3)$ implies
\be
Y(S^3, [g]) \geq Y(\RP^3, [g_0]) 
\ee
by (\ref{volYamabe}), which is a contradiction. Hence, there are no horizons. \stop

\section{An example with horizons} 

In this section, we provide an example to show that there exist metrics on $S^3$ with positive Ricci curvature such that the blowing-up manifolds  do have horizons.

Our example comes from a $1$-parameter family of left-invariant
metrics $\{ g_\ep \}$ on $S^3$, commonly known as the {\em Berger metrics}. Precisely, we think $S^3$ as the Lie Group
$$ SU(2) = 
\left\{
\left(
\begin{array}{cc}
z & - w\\
\bar{w} & \bar{z}
\end{array}
\right) \ : \ 
|z|^2 + |w|^2 = 1
\right\} ,
$$
where the Lie algebra of $SU(2)$ is spanned by
$$X_1 = 
\left(
\begin{array}{cc}
i & 0 \\
0 & - i
\end{array}
\right), \ 
X_2=
\left(
\begin{array}{cc}
0 & 1 \\
-1 & 0
\end{array}
\right), \ \mathrm{and} \  
X_3 =
\left(
\begin{array}{cc}
0 &  i \\
i  &  0
\end{array}
\right) .$$
Then $\{ g_\ep \}$ is defined by declaring $X_1, X_2, X_3$ to be orthogonal, $X_1$ to have length $\ep$ and $X_2, X_3$ to be unit vectors.  Note that scalar multiplication on $S^3 \subset \mathbb{C}^2$ corresponds to multiplication on the
left by  matrices 
$ \left( \begin{array}{cc}
e^{i \theta}& 0 \\
0  & e^{- i \theta}
\end{array} \right)$ 
on $SU(2)$, hence $X_1$ is exactly tangent to 
the circle fiber of the {\em Hopf fibration} 
$$ \pi: S^3 \longrightarrow S^2 = S^3/S^1  $$
and $g_\ep$ shrinks the circle fiber as $\ep \rightarrow 0$.
 One fact of $g_\ep$ for small $\ep$ is that all sectional
curvature of $(S^3, g_\ep)$ lies in the interval $[\ep^2, 4 - 3 \ep^2]$
(see \cite{Peterson}), in particular $g_\ep$ has positive Ricci curvature.  

\begin{prop}
Let $P \in S^3 $ be a fixed point and $G_\ep$ be the Green's function
of the conformal Laplacian of $g_\ep$ at $P$.
Then $(S^3 \setminus \{ P \}, G_\ep^4 g_\ep)$
has a horizon  for  $\ep$ sufficiently small.
\end{prop}

\pf For each $\ep \in (0, 1]$, we consider the rescaled metric $\bg_\ep = \ep^{-2} g_\ep$
and the Green's function $\bG_\ep$ associated to $\bg_\ep$ at $P$. 
Then, with respect to $\bg_\ep$, $X_1$ becomes a unit vector and 
$X_2, X_3$ have large
length $\ep^{-1}$ as $\ep \rightarrow 0$. Let $U \subset S^3$ be a 
fixed neighborhood of $P$ such that $\pi |_U$ is a trivial fiberation.
Let $O$ be a fixed point in the product manifold $S^1 \times \R^2$. 
By a scaling argument, there exists a family of diffeomorphisms  
$$
 \Psi_\ep :  U  \longrightarrow  \Psi_\ep(U) \subset S^1 \times \R^2 , 
$$
such that $\Psi_\ep (P) = O \in \Psi_\ep(U)$, $\{ \Psi_\ep(U) \}_{
1 \geq  \ep > 0} $ forms 
an exhaustion family of $S^1 \times \R^2$ as $\ep \rightarrow 0$, and
the push forward metrics 
$ \hg_\ep = \Psi_\ep^{-1 *}(\bg_\ep |_U)$ on $\Psi_\ep(U)$ 
converge in $C^2$ norm on compact sets to a flat metric $\hg$ on 
$S^1 \times \R^2$. Now fix another point  $Q \in \Psi_1(U)$ that is  
different from $O$ and consider the normalized function
\be
\hG_\ep (x) =  
\frac{ \bG_\ep \circ \Psi_\ep^{-1} (x) }{\bG_\ep \circ \Psi_\ep^{-1}(Q) }
\ee
for  $ x \in \Psi_\ep (U) \setminus \{O\}$. Then $\hG_\ep$ 
satisfies
\be
\left\{
\begin{array}{rcl}
8 \triangle_{\hg_\ep} \hG_\ep - R(\hg_\ep) \hG_\ep &  =  & 0 \ \  
\mathrm{on} \ \Psi_\ep(U) \setminus \{ O \} \\
\hG_\ep  & = & 1 \ \ \mathrm{at} \ Q 
\end{array} .
\right.
\ee
Since $\hG_\ep$ is positive and $\hg_\ep$ converges to 
$\hg$ as $\ep \rightarrow 0$, it follows from the Harnack inequality 
that $\hG_\ep$ coverges
to a positive function $\hG$ on $(S^1 \times \R^2) \setminus \{O\}$ in
$C^2$ norm on any compact set away from $\{O\}$. Furthermore, $\hG$ satisfies
\be
\left\{
\begin{array}{rcl}
 \triangle_{\hg } \hG  &  =  & 0 \ \  
\mathrm{on} \ (S^1 \times \R^2) \setminus \{ O \} \\
\hG & = & 1 \ \ \mathrm{at} \ Q
\end{array} .
\right.
\ee
On the other hand, the fact that the geodesic ball in  
$(S^1 \times \R^2, \hg)$ only has quadratic volume growth implies 
$(S^1 \times \R^2 , \hg)$ does not have a positive Green's function
for the usual Lapacian  $\triangle_{\hg}$ \cite{Cheng-Yau}. 
Therefore, $\hG \equiv 1$ on  $(S^1 \times \R^2) \setminus \{ O \}$.
Hence, the metrics $\hG_\ep^4 \hg_\ep $ converge to $\hg$ 
in $C^2$ norm on any compact set away from $\{O\}$. Now let $V
\subset S^1 \times \R^2$ be a small open ball containning $O$ such that $\partial V$ is an embedded two sphere whose
mean curvature vector computed with respect to $\hg$ points
towards $O$. Then, for suffiently small $\ep$, the mean curvature vector
of $\partial V$ computed with respect to $\hG^4_\ep \hg_\ep$ still points towards $O$. As $(\Psi_\ep (U), \hG^4_\ep \hg_\ep)$ is isometric to $(U, \bG^4_\ep \bg_\ep)$, the mean curvature vector of the boundary of  $\Psi_\ep^{-1}(V)$ in 
$(S^3 \setminus \{ P \}, \bG^4_\ep \bg_\ep)$ must  point towards
the blowing-up point $P$. On the other hand, as $(S^3 \setminus
\{ P \}, \bG^4_\ep \bg_\ep)$ is asymptotically flat,  its infinity 
is foliated by two shperes whose mean curvature vector points away from $P$. Therefore, it follows from standard geometric measure theory that there exists an embedded minimal two sphere in $\Psi_\ep (V)$, hence $(S^3 \setminus \{ P \}, \bG_\ep^4 \bg)$ has a horizon. \stop

\vh

\noindent {\bf Acknowledgment} 
I want to thank Justin Corvino and Rick Schoen for helpful discussions.

\bibliographystyle{plain}
\bibliography{RicMS}

%\vh
%Department of Mathematics, University of California, Santa Barbara,
%CA 93106
%{\em E-mail address: pengzim@math.ucsb.edu}

\end{document}